\documentclass{article}
\usepackage{latexsym,amsmath,amssymb,amscd}
\begin{document}

\title{On spherically symmetric solutions of the Einstein-Euler-de Sitter
equations   }
\author{Tetu Makino \footnote{Professor Emeritus at Yamaguchi University, Japan /e-mail: makino@yamaguchi-u.ac.jp}}
\date{\today}
\maketitle

\begin{abstract}
We construct spherically symmetric solutions to the Einstein-Euler equations, which 
contains a positive cosmological constant, say,
the Einstein-Euler-de Sitter equations. We assume a realistic barotropic equation of state. Equilibria of the spherically symmetric Einstein-Euler-de Sitter equations are given by the Tolman-Oppenheimer-Volkoff-de Sitter equation.
We can construct solutions near time periodic linearized solutions around the equilibrium. 
The Cauchy problem around the equilibrium can be solved.
This work can be considered as a trial of the generalization of the previous work on the
problem without cosmological constants.



\end{abstract}

\newtheorem{Lemma}{Lemma}
\newtheorem{Proposition}{Proposition}
\newtheorem{Theorem}{Theorem}
\newtheorem{Definition}{Definition}
\newtheorem{Assumption}{Assumption}
\newtheorem{Notation}{Notation}

\numberwithin{equation}{section}

\section{Introduction}

We consider the Einstein-de Sitter equation
$$R_{\mu\nu}-\frac{1}{2}g_{\mu\nu}(g^{\alpha\beta}R_{\alpha\beta})-\Lambda g_{\mu\nu}=
\frac{8\pi G}{c^4}T_{\mu\nu}
$$
for the energy-momentum tensor of a perfect fluid
$$T^{\mu\nu}=(c^2\rho +P)U^{\mu}U^{\nu}-Pg^{\mu\nu}.
$$
Here $R_{\mu\nu}$ is the Ricci tensor associated with the 
metric
$$
ds^2=g_{\mu\nu}dx^{\mu}dx^{\nu},$$
$G$ is the gravitational constant, $c$ the speed of light,
$\rho$ the mass density, $P$  the pressure and $U^{\mu}$ is 
the 4-dimensional velocity. $\Lambda$ is the cosmological
constant which is supposed to be positive in this article.  
See \cite[\S 111]{LandauL}.

Spherically symmetric solutions for the problem with
$\Lambda=0$ was investigated in \cite{ssEE}and \cite{JDE2017},  and the aim of
this article is to describe the similar results on
the problem with positive cosmological constants.\\

In this article we suppose that the pressure $P$ is a given function of
the density $\rho$ and suppose the following\\

\begin{Assumption}\label{Ass.1}
 $P$ is a smooth function of $\rho>0$ such that
$0<P, 0<dP/d\rho<c^2$ for $\rho>0$, and $P\rightarrow 0$ as
$\rho \rightarrow +0$. Moreover
there are positive constants $A, \gamma$ and a smooth function $\Omega$ defined on
$\mathbb{R}$ which is analytic on a neighborhood of $0$ such that 
$\Omega(0)=1$ and
$$P=A\rho^{\gamma}\Omega(A\rho^{\gamma-1}/c^2).
$$ 
for $\rho >0 $. 
We assume that $1<\gamma<2$.
\end{Assumption}

Moreover we suppose that either following Assumption 2 or Assumption 3 holds:

\begin{Assumption}\label{Ass.2}
$\displaystyle \frac{1}{\gamma-1}$ is an integer.
\end{Assumption}
\begin{Assumption}\label{Ass.3}
It holds that
$ \displaystyle 1<\gamma <\frac{54}{53}$.
\end{Assumption}

We consider spherically symmetric metrics of the form
$$ds^2=e^{2F(t,r)}c^2dt^2-
e^{2H(t,r)}dr^2-R(t,r)^2(d\theta^2+\sin^2\theta d\phi^2).
$$
We suppose that the system of coordinates is co-moving, that is,
$$U^0=e^{-F},\qquad U^1=U^2=U^3=0
$$
for $x^0=ct, x^1=r, x^2=\theta, x^3=\phi$. Then
the equations turn out to be
\begin{subequations}
\begin{eqnarray}
e^{-F}\frac{\partial R}{\partial t}&=&V \label{Aa} \\
e^{-F}\frac{\partial \rho}{\partial t}&=&
-(\rho+P/c^2)\Big(\frac{V'}{R'}+\frac{2V}{R}\Big) \label{Ab} \\
e^{-F}\frac{\partial V}{\partial t}&=&
-GR\Big(\frac{m}{R^3}+\frac{4\pi P}{c^2}\Big)+\frac{c^2\Lambda}{3}R+ \nonumber \\
&&-\Big(1+\frac{V^2}{c^2}-\frac{2Gm}{c^2R}
-\frac{\Lambda}{3}R^2\Big)\frac{P'}{R'(\rho+P/c^2)} \label{Ac} \\
e^{-F}\frac{\partial m}{\partial t}&=&
-\frac{4\pi}{c^2}R^2PV \label{Ad}
\end{eqnarray}
\end{subequations}
Here $X'$ stands for $\partial X/\partial r$.

The coefficients of the metric are given by
$$P'+F'(c^2\rho+P)=0$$
and 
$$e^{2H}=\Big(
1+\frac{V^2}{c^2}-\frac{2Gm}{c^2R}-\frac{\Lambda}{3}R^2\Big)^{-1}
(R')^2.
$$
In order to specify the function $F$, we introduce the state variable 
$u$ by
$$u=\int_0^{\rho}\frac{dP}{\rho+P/c^2}.
$$
Then there should exist a positive constant $C$ such that 
$$e^F=Ce^{-u/c^2}.
$$
To fix the idea we shall take $C=\sqrt{\kappa_+}$ with a positive constant $\kappa_+$ specified in the next Section. But the choice of the constant $C$ is free by change of the scale of $ct$.

We note that there are smooth functions $\Omega_u, 
\Omega_{\rho}, \Omega_P$ analytic on a neighborhood of
$0$ such that 
$\Omega_u(0)=\Omega_{\rho}(0)=\Omega_P(0)=1$ and
\begin{align*}
u&=\frac{\gamma A}{\gamma-1}\rho^{\gamma-1}\Omega_u(A\rho^{\gamma-1}/c^2), \\
\rho&=A_1u^{\frac{1}{\gamma-1}}\Omega_{\rho}(u/c^2), \\
P&=AA_1^{\gamma}u^{\frac{\gamma}{\gamma-1}}\Omega_P(u/c^2).
\end{align*}
Here $\displaystyle A_1:=\Big(\frac{\gamma-1}{\gamma A}\Big)^{\frac{1}{\gamma-1}}$. See \cite{TOVdS}.

We put
$$
m=4\pi \int_0^r\rho R^2R'dr,$$
supposing that $\rho$ is continuous at $r=0$. The coordinate $r$ can be 
changed to $m$, supposing that $\rho>0$, and the equations are reduced to
\begin{subequations}
\begin{eqnarray}
e^{-F}\Big(\frac{\partial R}{\partial t}\Big)_m&=&
\Big(1+\frac{P}{c^2\rho}\Big)V, \label{Ba} \\
e^{-F}\Big(\frac{\partial V}{\partial t}\Big)_m&=&
\frac{4\pi}{c^2}R^2PV\frac{\partial V}{\partial m}-GR\Big(\frac{m}{R^3}+\frac{4\pi P}{c^2}\Big)+\frac{c^2\Lambda}{3}R+ \nonumber \\
&&-\Big(1+\frac{V^2}{c^2}-\frac{2Gm}{c^2R}-\frac{\Lambda}{3}R^2\Big)
\Big(1+\frac{P}{c^2\rho}\Big)^{-1}
4\pi R^2\frac{\partial R}{\partial m}. \label{Bb}
\end{eqnarray}
\end{subequations}
Here $(\partial/\partial t)_m$ means the differentiation with respect to $t$
keeping $m$ constant. We will change the coordinate $m$ to $r$ later
through a fixed equilibrium, and we shall construct solutions near the equilibrium.

\section{Equilibrium}

Let us consider solutions independent of $t$, that is,
$F=F(r), H=H(r), \rho=\rho(r), V\equiv 0, R\equiv r$. The
equations are reduced to the Tolman-Oppenheimer-Volkoff-de
Sitter equation
\begin{subequations}
\begin{eqnarray}
\frac{dm}{dr}&=&4\pi r^2\rho, \label{Ca} \\
\frac{dP}{dr}&=&-(\rho+P/c^2)
\frac{\displaystyle G\Big(m+\frac{4\pi r^3}{c^2}P\Big)-\frac{c^2\Lambda}{3}r^3}{r^2\displaystyle\Big(1-\frac{2Gm}{c^2r}-\frac{\Lambda}{3}r^2\Big)}. \label{Cb}
\end{eqnarray}
\end{subequations}
This equation was analyzed in \cite{TOVdS}. Let us summarize the results.

For arbitrary positive central density $\rho_{\mathsf{O}}$ there
exists a unique solution germ $(m(r), P(r)), 0<r\ll 1,$ such that
\begin{subequations}
\begin{eqnarray}
m&=&\frac{4\pi}{3}\rho_{\mathsf{O}}r^3(1+O(r^2)), \label{mc} \\
P&=&P_{\mathsf{O}}-(\rho_{\mathsf{O}}+P_{\mathsf{O}}/c^2)
\Big(\frac{4\pi G}{3}(\rho_{\mathsf{O}}+3P_{\mathsf{O}}/c^2)-\frac{c^2\Lambda}{3}\Big)
\frac{r^2}{2}+O(r^4). \label{Pc}
\end{eqnarray}
\end{subequations}

We denote
\begin{align}
\kappa(r.m)&:=1-\frac{2Gm}{c^2r}-\frac{\Lambda}{3}r^2, \\
Q(r, m, P)&:=G\Big(m+\frac{4\pi r^3}{c^2}P\Big)-\frac{c^2\Lambda}{3}r^3.
\end{align}

We concentrate ourselves to solutions satisfying $\kappa(r, m(r))>0$.

Let us use the concept of `monotone short solutions' defined by the following definition as \cite{TOVdS}:

\begin{Definition}
A a solution $(m(r), P(r)), 0<r<r_+,$ of (\ref{Ca})(\ref{Cb}) is
said to be {\bf monotone-short} if $r_+<\infty$,
$dP/dr<0$ for $0<r<r_+$, that is, $Q(r, m(r), P(r))>0$, and
$P \rightarrow 0$ as $r\rightarrow r_+-0$ and if both
\begin{equation}
\kappa_+:=\lim_{r\rightarrow r_+-0}\kappa(r,m(r))=
1-\frac{2Gm_+}{c^2r_+}-\frac{\Lambda}{3}r_+^2
\end{equation}
and 
\begin{equation}
Q_+:=\lim_{r\rightarrow r_+-0}Q(r, m(r), P(r))=
Gm_+-\frac{c^2\Lambda}{3}r_+^3
\end{equation}
are positive. Here 
\begin{equation}
m_+:=\lim_{r\rightarrow r_+-0}m(r)=4\pi\int_0^{r_+}
\rho(r)r^2dr.
\end{equation}
\end{Definition}

So, we suppose

\begin{Assumption}\label{Ass.4}
There is a monotone-short solution $(\bar{m}(r), \bar{P}(r)), 0<r<r_+,$ satisfying (\ref{mc})(\ref{Pc}).
\end{Assumption}

We fix the equilibrium given by the monotone-short solution hereafter.

As for sufficient conditions for the existence of
monotone-short prolongations, see \cite{TOVdS}.
The sufficient conditions given by \cite[Theorem 1]{TOVdS} and \cite[Theorem 2]{TOVdS} assume relative smallness of $\Lambda$. The task to prove the existence of 
monotone-short solutions seems to be not so easy when $\Lambda$ is permitted to be large. Actually, when $\Lambda=0$, $P(r)$ given by shooting from $r=0$ remains monotone decreasing with respect to $r$ automatically, and $\kappa_+$ turns out to be positive automatically if $P(r)$ hits zero at $r_+ <\infty$. But generally it is not the case when $\Lambda >0$ is large.

Anyway, the associated function $u=\bar{u}(r)$ turns out to be of the form
\begin{equation}
\bar{u}(r)=\frac{Q_+}{r_+^2\kappa_+}(r_+-r)(1+[r_+-r,
(r_+-r)^{\frac{\gamma}{\gamma-1}}]_1) \label{ur+}
\end{equation}
as $r\rightarrow r_+-0$. See \cite[Theorem 4]{TOVdS}. Here and hereafter we use

\begin{Notation}
The symbol $[X]_Q$ stands for various convergent power series of the form
$$\sum_{k\geq Q}a_kX^k.$$
The symbol $[X_1,X_2]_Q, Q=0,1,2,\cdots$ stands for various convergence double power series of the form
$$\sum_{k_1+k_2\geq Q}a_{k_1k_2}X_1^{k_1}X_2^{k_2}.$$
\end{Notation}

\section{Equations for the small perturbation from the equilibrium}

Using the fixed equilibrium $m=\bar{m}(r)$, we take the variable $r$
given by its inverse function.
We are going to a solution near equilibrium of the form
$$R=r(1+y),\qquad V=rv. $$
Here $y, v$ are small unknowns. The equations turn out to be
\begin{subequations}
\begin{eqnarray}
e^{-F}\frac{\partial y}{\partial t}&=&
\Big(1+\frac{P}{c^2\rho}\Big)v, \label{Ea} \\
e^{-F}\frac{\partial v}{\partial t}&=&
\frac{(1+y)^2}{c^2}\frac{P}{\bar{\rho}}v\frac{\partial}{\partial r}(rv) + \nonumber \\
&&-G(1+y)\Big(\frac{m}{r^3(1+y)^3}+\frac{4\pi}{c^2}P\Big)+
\frac{c^2\Lambda}{3}(1+y)+ \nonumber \\
&&-\Big(1+\frac{r^2v^2}{c^2}-
\frac{2Gm}{c^2r(1+y)}-
\frac{\Lambda}{3}r^2(1+y)^2\Big)
\times \nonumber \\
&&\times
\Big(1+\frac{P}{c^2\rho}\Big)^{-1}
\frac{(1+y)^2}{\bar{\rho}r}\frac{\partial P}{\partial r}. \label{Eb}
\end{eqnarray}
\end{subequations}

Here $m=\bar{m}(r)$ is a given function and $\rho, P$ are considered
as given functions of $r$ and the unknowns $y, z(:=r\partial y/\partial r)$ 
as follows:
\begin{align*}
\rho&=\bar{\rho}(r)(1+y)^{-2}(1+y+z)^{-1}, \\
P&=\bar{P}(r)(1-\Gamma (\bar{u}(r))(3y+z)-\Phi(\bar{u}(r), y, z)).
\end{align*}
Here
$$\Gamma:=\frac{\rho}{P}\frac{dP}{d\rho}
$$
and $\Phi(u, y, z)$ is an analytic function of the form
$\sum_{k_0\geq 0, k_1+k_2\geq 2}u^{k_0}y^{k_1}z^{k_2}$.
We shall denote such a function by $[u;y,z]_{0;2}$ hereafter.
Moreover we shall use
\begin{align*}
1+\frac{P}{c^2\rho}&=
\overline{\Big(1+\frac{P}{c^2\rho}\Big)}
\Big(1-\frac{\bar{P}}{c^2\bar{\rho}}
\overline{\Big(1+\frac{P}{c^2\rho}\Big)^{-1}}(\Gamma-1)(3y+z)+ \\
&+[\bar{u}(r); y, z]_{0;2}\Big).
\end{align*}

\section{Analysis of the linearized equation}

Let us linearize (\ref{Ea})(\ref{Eb}):
\begin{subequations}
\begin{eqnarray}
e^{-\bar{F}}\frac{\partial y}{\partial t}&=&
\overline{\Big(1+\frac{P}{c^2\rho}\Big)}v, \\
e^{-\bar{F}}\frac{\partial v}{\partial t}&=&E_2y''+E_1y'+E_0y,
\end{eqnarray}
\end{subequations}
where $y''=\partial^2y/\partial r^2, y'=\partial y/\partial r$ and
\begin{align*}
E_2&=e^{-2\bar{H}}\overline{(\rho+P/c^2)^{-1}}\overline{P\Gamma}, \\
\frac{E_1}{E2}&=
\frac{d}{dr}\Big(\bar{H}+\bar{F}-\log(\overline{1+P/c^2\rho})+
\log(\overline{P\Gamma}r^4)\Big), \\
E_0&= \frac{4\pi G}{c^2}3\overline{(\Gamma-1)P}+ \\
&+\Big(-1-3\overline{\Gamma e^{-2H}}+
3\overline{(\Gamma-1)e^{-2H}(1+P/c^2\rho)^{-1}}\Big)
 \overline{(\rho+P/c^2)^{-1}}\frac{1}{r}\frac{d\bar{P}}{dr} +\\
&+3e^{-2\bar{H}}\overline{(\rho+P/c^2)^{-1}}
\frac{1}{r}\frac{d}{dr}\overline{P\Gamma}+ \\
&+\Lambda\Big(c^2+r\frac{d\bar{u}}{dr}\Big).
\end{align*}
Here $\bar{X}, \overline{XXX}$ denote the evaluations along the fixed equilibrium. Putting
$$\mathcal{L}(D)y:=-e^{2\bar{F}}\overline{(1+P/c^2\rho)}
(E_2D^2y'+E_1Dy+E_0y),
$$
we get the linearized wave equation
\begin{equation}
\frac{\partial^2y}{\partial t^2}+\mathcal{L}\Big(\frac{\partial}{\partial r}\Big)y=0.
\end{equation}
We can rewrite $\mathcal{L}$ in the formally self-adjoint form
$$\mathcal{L}\Big(\frac{d}{dr}\Big)y=-\frac{1}{b}\frac{d}{dr}a\frac{dy}{dr}+Qy,
$$
where
\begin{subequations}
\begin{align}
a&=e^{\bar{H}+\bar{F}}\frac{\overline{P\Gamma}r^4}{\overline{1+P/c^2\rho}}
\\
b&=e^{3\bar{H}-\bar{F}}
\frac{r^4\bar{\rho}}{\overline{1+P/c^2\rho}} \\
Q&=-e^{2\bar{F}}\overline{1+P/c^2\rho }E_0.
\end{align}
\end{subequations}
It is easy to see that $Q$ is bounded on $0\leq r\leq r_+$.\\

For the sake of the simplicity of discussion, we suppose

\begin{Assumption}
The function $\rho \mapsto P$ is analytic near $\rho_{\mathsf{O}}$.
\end{Assumption}

As in \cite{ssEE}, we introduce the new variable $x$ instead of $r$ defined by
\begin{equation}
x:=\frac{\tan^2\theta}{1+\tan^2\theta} \quad\mbox{with}\quad
\theta:=\frac{\pi}{2\xi_+}\int_0^r
\sqrt{\frac{\bar{\rho}}{\overline{\Gamma P}}}
e^{-\bar{F}+\bar{H}}dr
\end{equation}
Here
\begin{equation}
\xi_+:=\int_0^{r_+}
\sqrt{\frac{\bar{\rho}}{\overline{\Gamma P}}}
e^{-\bar{F}+\bar{H}}dr.
\end{equation}

We can and shall assume that $\xi_+/\pi =1$ without loss of generality, by changing the scale of $ct$ or the coefficient of $e^{-\bar{F}}$.

Then there are positive constants $C_0, C_1$ such that
\begin{align}
r&=C_0\sqrt{x}(1+[x]_1) \quad \mbox{for}\quad 0<x\ll 1, \\
r_+-r&=C_1(1-x)(1+[1-x, [1-x]^{\frac{\gamma}{\gamma-1}}]_1)\quad
\mbox{for} \quad 0<1-x \ll 1
\end{align}
Using this variable, we can write the operator $\mathcal{L}$ as
\begin{align}
\mathcal{L}y&=
-x(1-x)\frac{d^2y}{dx^2}-
\Big(\frac{5}{2}(1-x)-\frac{N}{2}x\Big)\frac{dy}{dx}+ \nonumber \\
& +L_1(x)\frac{dy}{dx}+L_0(x)y, \label{L}
\end{align}
where $L_1(x), L_0(x)$ are smooth functions on $[0,1[$ such that
$$L_{\mu}(x)=
\begin{cases}
[x]_{\mu}\quad\mbox{for}\quad 0<x \ll 1 \\
[1-x,(1-x)^{\frac{N}{2}}]_{\mu}\quad\mbox{for}\quad 0<1-x\ll 1
\end{cases}
$$
for $\mu=0,1$. Here and hereafter we use
\begin{Notation}
Put
$$N=\frac{2\gamma}{\gamma-1}.$$
\end{Notation}

Assumption \ref{Ass.2} means that $N$ is an even integer and Assumption \ref{Ass.3}
means that $N >108$.

Therefore \cite[Proposition 7]{ssEE} is still valid, that is, we can claim

\begin{Proposition}
The operator $\mathfrak{T}_0,
\mathcal{D}(\mathfrak{T}_0)=
C_0^{\infty}(]0,r_+[)$, $\mathfrak{T}_0y=\mathcal{L}y$ in the Hilbert space
$L^2([0,r_+], b(r)dr)$ admits the Friedrichs extension $\mathfrak{T}$, a
self adjoint operator, whose spectrum consists of
simple eigenvalues $\lambda_1<\lambda_2<\cdots<\lambda_{\nu}<\cdots
\rightarrow +\infty$.
\end{Proposition}

 Since the discussion is
quite parallel to that of \cite[\S 5]{ssEE} and \cite{JDE2017}, we omit the details.

Hereafter we shall denote by $\mathcal{L}$ the self-adjoint operator $\mathfrak{T}$ in
$L^2([0,r_+], b(r)dr)$.
Moreover it is easy to see, keeping in mind the behavior of $\bar{u}(r)$ as $r\rightarrow r_+$ given by \eqref{ur+},  that $L^2([0,r_+], b(r)dr)$ is isomorphic to 
$L^2([0,1], x^{\frac{3}{2}}(1-x)^{\frac{N}{2}-1}dx)$, 
the space of functions of the variable $x$. In this sense we consider $\mathcal{L}$ as a self-adjoint operator acting on functions of $x$ in 
$L^2([0,1], x^{\frac{3}{2}}(1-x)^{\frac{N}{2}-1}dx)$ of the form \eqref{L}.\\

As \cite[Proposition 10]{ssEE} and \cite[Appendix A]{JDE2017}, we can claim 

\begin{Proposition}
Any eigenfunction $\psi(x)$ of $\mathcal{L}$ is of the form
$$
\psi(x)=
\begin{cases}
C_0 +[x]_1 \quad \mbox{for}\quad 0<x\ll 1 \\
C_1+[1-x, (1-x)^{N/2}]_1
\quad\mbox{for}\quad 0< 1-x \ll 1
\end{cases},
$$
where the constants $C_0, C_1$ are such that $C_0\not=0, C_1\not=0$.
\end{Proposition}

Proof is due to \cite[Lemma 2]{OJM}.

\section{Rewriting the equation system (\ref{Ea})(\ref{Eb}) using $\mathcal{L}$}

Let us go back to the system of equations
(\ref{Ea})(\ref{Eb}). 
We are going to rewrite the system of equations \eqref{Ea} \eqref{Eb} as
\begin{subequations}
\begin{align}
&\frac{\partial y}{\partial t}-
J(x,y,x\frac{\partial y}{\partial x})v=0, \label{EVa} \\
&\frac{\partial v}{\partial t}+H_1
(x,y,x\frac{\partial y}{\partial x}, v, \Lambda)\mathcal{L}\Big(\frac{\partial}{\partial x}\Big)
y+H_2(x,y,x\frac{\partial y}{\partial x}, v, x\frac{\partial v}{\partial x}, \Lambda)=0 \label{EVb}
\end{align}
\end{subequations}\\

We shall use the analysis of
$\partial P/\partial r$ given in \cite[(6.2)]{ssEE}:
\begin{align*}
-\frac{1}{r\bar{\rho}}\frac{\partial P}{\partial r}&=
-\frac{1}{r\bar{\rho}}\frac{d\bar{P}}{dr}+
(1+\partial_z\Phi/\Gamma)\frac{1}{r\bar{\rho}}
\frac{\partial}{\partial r}\overline{P\Gamma}(3y+z) + \\
&+\frac{\bar{P}}{r\bar{\rho}}\cdot[Q0]+
\frac{1}{r\bar{\rho}}\frac{d\bar{P}}{dr}\cdot
[Q1],
\end{align*}
where $[Q0], [Q1]$ are given by \cite[(6.3a)]{ssEE}, \cite[(6.3b)]{ssEE}.

(The symbols $\Phi^P, \gamma^P$ of
\cite{ssEE} are replaced by $\Phi, \Gamma$ in this article.)

We put
$$\mbox{the right-hand side of (\ref{Eb})}=[R2]+[R1]+[W]\frac{\Lambda}{3},
$$
where
\begin{align*}
[R2]&:=\frac{(1+y)^2}{c^2}\frac{P}{\bar{\rho}}v(v+w)\quad\mbox{with}\quad w=r\frac{\partial v}{\partial r}, \\
[W]&:=c^2(1+y)-
r^2(1+y)^4(1+P/c^2\rho)^{-1}\Big(-\frac{1}{r\bar{\rho}}\frac{\partial P}{\partial r}\Big).
\end{align*}
We put
$$ [R1]=[R3]+[R4]+[R5]+[R6]+[R7] $$
as in \cite{ssEE}.  But the analysis of $[W]$ is new:
We put
\begin{align*}
[W]&=c^2+r\frac{d\bar{u}}{dr}+[W1]+[W2]+[W3]+[W4], \\
[W1]&:=c^2y-r^2(1+y)^2(1+P/c^2\rho)^{-1}
\Big(-\frac{1}{r\bar{\rho}}\frac{d\bar{P}}{dr}\Big)-r\frac{d\bar{u}}{dr} \\
&=[W1L]+[W1Q], \\
[W1L]&:=c^2y-4r\frac{d\bar{u}}{dr}y+
r\overline{\Big(\frac{P/c^2}{\rho+P/c^2}\Big)}
(\Gamma-1)(3y+z)\frac{d\bar{u}}{dr}, \\
[W2]&:=-r^2(1+y)^4(1+P/c^2\rho)^{-1}
(1+\partial_z\Phi/\Gamma)\frac{1}{r\bar{\rho}}
\frac{\partial}{\partial r}\overline{P\Gamma}(3y+z), \\
[W3]&:=-r^2(1+y)^4(1+P/c^2\rho)^{-1}\frac{\bar{P}}{r\bar{\rho}}
[Q0], \\
[W4]&:=-r^2(1+y)^4(1+P/c^2\rho)^{-1}
\frac{1}{r\bar{\rho}}
\frac{d\bar{P}}{dr}[Q1].
\end{align*}

Then it follows from (\ref{Cb}) that
\begin{align*}
[R1]+[W]\frac{\Lambda}{3}&=[R3L]+[R3Q]+[R4L]+[R4Q]+ [R5]+[R6]+[R7]+ \\
&+([W1L]+[W1Q]+[W2]+[W3]+[W4])\frac{\Lambda}{3}.
\end{align*}

Let us define
\begin{align*}
1+G_1&=(1+\partial_z\Phi/\Gamma)\Big(1+\frac{r^2v^2}{c^2}-
\frac{2Gm}{c^2r(1+y)}-r^2(1+y)^2\frac{\Lambda}{3}\Big)\times \\
&\times\Big(1-\frac{2Gm}{c^2r}-\frac{\Lambda}{3}r^3\Big)^{-1}
\frac{1+\overline{P/c^2\rho}}{1+P/c^2\rho}
(1+y)^2.
\end{align*}
Then we have
\begin{align*}
-e^{-2\bar{F}}\overline{(1+P/c^2\rho)^{-1}}\mathcal{L}y&=
[R3L]+[R4L]+[W1L]\frac{\Lambda}{3} + \\
&+\frac{1}{1+G_1}\Big([R5]+[W2]\frac{\Lambda}{3}\Big).
\end{align*}
Putting
\begin{align*}
G_2&:=G_1([R3L]+[R4L]+[W1L]\frac{\Lambda}{3})+ \\
&-([R3Q]+[R4Q]+[R6]+[R7]+[R2]) + \\
&-([W1Q]+[W3]+[W4])\frac{\Lambda}{3}, \\
H_2&:=e^FG_2, \\
H_1&:=e^{F-2\bar{F}}
\overline{(1+P/c^2\rho)^{-1}}(1+G_1),
\end{align*}
we can write
$$e^F\times(\mbox{the right-hand side of (\ref{Eb})})=-H_1\mathcal{L}y-H_2.
$$

Putting
$$J:=e^F(1+P/c^2\rho),$$
we rewrite the system of equations (\ref{Ea})(\ref{Eb}) as
\eqref{EVa}\eqref{EVb}.\\

Then we see that  the functions
\begin{align*}
&J : (x,y,z) \mapsto J(x,y, z,\Lambda) \\
&H_1: (x,y,z,v) \mapsto H_1(x,y,z,v,\Lambda) \\
&H_2 : (x,y,z,v,w) \mapsto H_2(x,y,z,v,w,\Lambda),
\end{align*}
$\Lambda$ being fixed, enjoy the following properties:\\

{\bf (B.1)}: {\it $J, H_1, H_2$ belong to the functional classes
$\mathfrak{A}_{(N)}^0(U^2), \mathfrak{A}_{(N)}^0(U^3),\mathfrak{A}_{(N)}^2(U^4)$, respectively.} \\

Here we use

\begin{Definition}
By $\mathfrak{A}_{(N)}^Q(U^p)$, $U$ being the interval $]-\delta, \delta[, 0<\delta \ll 1$,
we denote the set of all smooth functions $f$ defined on $[0, 1[\times U^p$ such that there are convergent power series
$$
\Phi_0(X,Y_1,\cdots, Y_p)=\sum_{|\vec{k}|\geq Q}
a_{j\vec{k}}X^jY_1^{k_1}\cdots Y_p^{k_p} $$
and
$$
\Phi_1(X_1,X_2,Y_1,\cdots, Y_p)=
\sum_{|\vec{k}|\geq Q}b_{j_1j_2\vec{k}}
X_1^{j_1}X_2^{j_2}Y_1^{k_1}\cdots Y_p^{k_p}
$$
such that
$$f(x,y_1,\cdots, y_p)=\Phi_0(x,y_1,\cdots, y_p)\quad\mbox{for}\quad
0<x \ll 1$$
and
$$f(x,y_1, \cdots, y_p)=\Phi_1(1-x,(1-x)^{N/2},y_1,\cdots, y_p)\quad\mbox{for}\quad
0<1-x \ll 1.$$
\end{Definition}

Clearly  \\

{\bf (B.2)}:  {\it It holds that
\begin{equation}
J(x,0,0)H_1(x,0,0,0)=1, 
\end{equation}
and
\begin{equation}
\frac{1}{C}<J(x, 0, 0) <C 
\end{equation}
with a sufficiently large $C$}. \\

Moreover we see \\

{\bf (B.3)}: {\it It holds that }
\begin{subequations}
\begin{align}
&\frac{\partial J}{\partial z}\equiv_{(N)}0, \\
&\Big(\frac{\partial H_1}{\partial z}\Big)\mathcal{L}y+
\frac{\partial H_2}{\partial z}\equiv_{(N)}0, \label{HLH=0} \\
&\frac{\partial H_2}{\partial w}\equiv_{(N)}0.
\end{align}
\end{subequations}\\

Here we mean

\begin{Notation}
For $f \in \mathfrak{A}_{(N)}^0(U^p)$, `$f\equiv_{(N)}0$' means
that there exists a convergent power series 
$\Phi(X_1,X_2,Y_1,\cdots, Y_p)$ such that
$$f(x,y_1,\cdots, y_p)=(1-x)
\Phi(1-x,(1-x)^{N/2},y_1,\cdots, y_p)
\quad\mbox{for}\quad 0<1-x \ll 1.$$
\end{Notation}

Proof of \eqref{HLH=0}:
$$(\partial_zH_1)\mathcal{L}y+\partial_zH_2 \equiv_{(N)} 0  $$
 is similar to that of \cite[Proposition 11]{ssEE}. 

Actually we see
$$(\partial_zH_1)\mathcal{L}y+\partial_zH_2\equiv_{(N)} e^F[S]$$
and we have to show $[S]\equiv_{(N)} 0$, where
\begin{align*}
[S]&:=(\partial_zG_1)\Big(e^{-2\bar{F}}
\overline{(1+P/c^2\rho)^{-1}}
\mathcal{L}y+[R3L]+[R4L]+[W1L]\Lambda/3\Big)+ \\
&+G_1\partial_z([R3L]+[R4L]+[W1L]\Lambda/3)+ \\
&-\partial_z\Big([R3Q]+[R4Q]+[R6]+[R7]+[R2]+ ([W1Q]+[W3]+[W4])\Lambda/3\Big).
\end{align*}
But we have
$$[S]\equiv_{(N)}-\frac{\partial_zG_1}{1+G_1}\Big([R5]+[W2]\frac{\Lambda}{3}\Big)-
\partial_z\Big([R7]+[W4]\frac{\Lambda}{3}\Big),
$$
since $\partial_z[R3L]$, $\partial_z[R4L]$, $\partial_z[R3Q]$,
$\partial_z[R4Q]$, $\partial_z[R6]$, $\partial_z[R2]$,
$\partial_z[W1L]$, $\partial_z[W1Q]$, $\partial_z[W3]$ are all $\equiv_{(N)} 0$ clearly.
By a tedious calculation, we get
\begin{align*}
&-\frac{\partial_zG_1}{1+G_1}\Big([R5]+[W2]\frac{\Lambda}{3}\Big)\equiv
\partial_z\Big([R7]+[W4]\frac{\Lambda}{3}\Big) \\
&\equiv_{(N)}
-\partial_z^2\Phi\Big(1+\frac{r^2v^2}{c^2}
-\frac{2Gm}{c^2r(1+y)}-r^2(1+y)^2\frac{\Lambda}{3}\Big)(1+y)^2
\frac{1}{r\bar{\rho}}
\frac{d\bar{P}}{dr}(3y+z),
\end{align*}
so that $[S]\equiv_{(N)} 0$. This completes the proof.

\section{Main results}

Let us fix a time periodic solution of the
linearized equation:
$$Y_1=\sin(\sqrt{\lambda}t+\Theta_0)\psi(x),$$
where $\lambda$ is a positive eigenvalue of the operator
$\mathcal{L}$, $\psi$ is an associated eigenfunction
and $\Theta_0$ is an arbitrary constant. 
We seek a solution of the form
\begin{equation}
\begin{bmatrix}
y \\
\\
v
\end{bmatrix}
=\begin{bmatrix}
\varepsilon(Y_1+\check{y}) \\
\\
\varepsilon(V_1+\check{v})
\end{bmatrix}
=\varepsilon
\begin{bmatrix}
Y_1 \\
\\
V_1
\end{bmatrix}
+\varepsilon\vec{w}, \label{Defw}
\end{equation}
where $$V_1=e^{-\bar{F}}(1+\overline{P/c^2\rho})^{-1}
\frac{\partial Y_1}{\partial t}.$$
Then  the equation to be solved turns out to be 
\begin{equation}
\mathfrak{P}(\vec{w})=\varepsilon \vec{c}, \label{P}
\end{equation}
for the unknown $\vec{w}= (\check{y}, \check{v})^{\top}$.
Here 
$$\mathfrak{P}(\vec{w})-\varepsilon\vec{c}=
\frac{1}{\varepsilon}
\begin{bmatrix}
\mbox{the left-hand side of \eqref{EVa} } \\
\\
\mbox{the left-hand side of \eqref{EVb} }
\end{bmatrix}_{y=\varepsilon(Y_1+\check{y}), v=\varepsilon(V_1+\check{v})}
.
$$
We are seeking the solution $\vec{w}$ satisfying the initial condition
\begin{equation}
\vec{w}|_{t=0}=
\begin{bmatrix}
\check{y} \\
\\
\check{v}
\end{bmatrix}_{t=0}=\vec{0}. \label{Pinit}
\end{equation}\\

The Fr\'{e}chet derivative $D\mathfrak{P}(\vec{w})$ of the nonlinear operator
$\mathfrak{P}$ at $\vec{w}$ :
$$D\mathfrak{P}(\vec{w})\vec{h}=
\begin{bmatrix}
\underline{DP1} \\
\\
\underline{DP2}
\end{bmatrix},
\quad\mbox{with}\quad
\vec{h}=
\begin{bmatrix}
h \\
k
\end{bmatrix}
$$
is given by
\begin{align*}
\underline{DP1}&=\displaystyle\frac{\partial h}{\partial t}
-Jk-\Big(
(\partial_yJ)v+(\partial_zJ)vx\frac{\partial}{\partial x}\Big)h \\
\underline{DP2}&=\displaystyle\frac{\partial k}{\partial t}
+H_1\mathcal{L}h+ \\
&+\Big((\partial_yH_1)\mathcal{L}y+\partial_yH_2+
((\partial_zH_1)\mathcal{L}y+\partial_zH_2)x\frac{\partial}{\partial x}\Big)h + \\
&+\Big((\partial_vH_1)\mathcal{L}y+\partial_vH_2+
\partial_wH_2x\frac{\partial}{\partial x}\Big)k.
\end{align*}

Thanks to {\bf (B.3)} we can claim that 
there are  functions $a_{01}$, $a_{00}$,
$a_{11}$, $a_{10}$, $a_{21}$, $a_{20}$ of $t, x,y, \partial_xy, \partial_x^2y,
v, \partial_xv$ in the class $\mathfrak{A}_{(N)}([0,T]\times U^5)$ such that
\begin{align*}
\underline{DP1}&=\frac{\partial h}{\partial t}-Jk+\Big(a_{01}x(1-x)\frac{\partial}{\partial x}+a_{00}\Big)h, \\
\underline{DP2}&=\frac{\partial k}{\partial t}+H_1\mathcal{L}h+
\Big(a_{11}x(1-x)\frac{\partial}{\partial x}+a_{10}\Big)h + \Big(a_{21}x(1-x)\frac{\partial}{\partial x}+a_{20}\Big)k.
\end{align*}

Here we use 

\begin{Definition}
By $\mathfrak{A}_{(N)}([0,T]\times U^p)$ we denote the set of all smooth functions $a(t,x,y_1,\cdots, y_p)$ of 
$t \in [0,T], x \in [0,1[, y_j \in U=]-\delta,\delta[, 0<\delta \ll 1$, such that there are analytic functions $\Psi_0$ on $[0,T]\times U^{p+1}$, $\Psi_1$ on
$[0,T] \times U^{p+2}$ such that
$$a(t,x,y_1,\cdots)=
\begin{cases}
\Psi_0(t,x,y_1,\cdots)\quad\mbox{for}\quad 0<x \ll 1 \\
\Psi_1(t,1-x,(1-x)^{N/2}, y_1, \cdots)\quad\mbox{for}\quad 
0<1-x\ll 1
\end{cases}.
$$
\end{Definition}

In fact {\bf (B.3)} plays a crucial r\^{o}le. Otherwise the factor $(1-x)$ in terms
$a_{01}x(1-x)\frac{\partial}{\partial x}, a_{11}x(1-x)\frac{\partial}{\partial x}, a_{21}x(1-x)\frac{\partial}{\partial x} $ would lack so that these terms lacking the factor $1-x$ would invade the principal part of the linearized operator $\mathcal{L}$. 

Then we can apply the Nash-Moser theorem to get the main results,
since the Fr\'{e}chet derivative $D\mathfrak{P}(\vec{w})$ has the inverse, say, the resolution of the wave equation,  on a suitable graded functional spaces of functions $\vec{h}$ such that $\vec{h}|_{t=0}=0$. We employ the Nash-Moser theorem, since this resolution involves so called regularity loss at the vacuum boundary $x=1$, while
the singularity at the center $x=0$ does not cause regularity loss because of the factor
$x$ in terms $\displaystyle z=x\frac{\partial y}{\partial x}, w=x\frac{\partial v}{\partial x}$. 
 
Namely, 
under the Assumption \ref{Ass.2}, when $N/2$ is an integer, we apply the Nash-Moser theorem formulated by R. S. Hamilton given in \cite{Hamilton} as \cite{ssEE}, while,
under the Assumption \ref{Ass.3}, when $N >108$, we apply the Nash-Moser theorem formulated by J. T. Schwartz in \cite{Schwartz} as \cite{JDE2017}. 

The results are following:

\begin{Theorem}
Given $T>0$ under Assumption \ref{Ass.2}, there is a positive number $\epsilon_0$ such that,
for
$|\varepsilon|\leq \epsilon_0$, there is a solution
$\vec{w}\in C^{\infty}([0,T]\times [0,1])$ of \eqref{P}\eqref{Pinit} such that
$$\sup_{j+k\leq n}
\Big\|\Big(\frac{\partial}{\partial t}\Big)^j\Big(\frac{\partial}{\partial x}\Big)^k\vec{w}\Big\|_{L^{\infty}([0,T]\times[0,1])}
\leq C(n)|\varepsilon|, \qquad \forall n \in \mathbb{N},$$
and hence a solution $(y,v)$ of (\ref{Ea})(\ref{Eb}) of the form
$y=\varepsilon Y_1+O(\varepsilon^2),
v=\varepsilon V_1+O(\varepsilon^2)$.
\end{Theorem}

\begin{Theorem}
Given $T>0$ under Assumption \ref{Ass.3}, there is a positive number $\epsilon_0$ such that,
for
$|\varepsilon|\leq \epsilon_0$, there is a solution
$\vec{w}$ of \eqref{P}\eqref{Pinit} such that
$$\|\vec{w}\|_{s_N+1} \leq C\varepsilon,$$
where $\displaystyle s_N=\Big[\frac{N}{2}\Big]+1=\mathrm{min}\Big\{ s \in \mathbb{N}
\Big|
s>\frac{N}{2}\Big\}$. 
\end{Theorem}

Here we use

\begin{Definition}
We put
$$\|\vec{w}\|_{\nu}^2=\|y\|_{\nu}^2+\|v\|_{\nu}^2$$
for $\vec{w}=(y,v)^{\top}$, where
\begin{align*}
&\|u\|_{\nu}^2=\sum_{\iota+\kappa \leq \nu}
\int_0^T(\|(-\partial_t)^{2\iota}u(t,\cdot)\|_{\kappa}^*)^2dt, \\
&(\|u\|_{\kappa}^*)^2=(\|u^{[0]}\|_{[0]\kappa}^*)^2+
(\|u^{[1]}\|_{[1]\kappa}^*)^2, \\
&(\|u^{[\mu]}\|_{[\mu]\kappa}^*)^2=\sum_{0\leq m\leq \kappa}
\|\triangle_{[\mu]}^mu^{[\mu]}\|_{[\mu]}^2, \quad \mu=0,1\\
&\triangle_{[0]}=x\frac{d^2}{dx^2}+\frac{5}{2}\frac{d}{dx}, \qquad
\triangle_{[1]}=X\frac{d^2}{dX^2}+\frac{N}{2}\frac{d}{dX}\quad\mbox{for}\quad X=1-x, \\
&\|f\|_{[0]}^2=\int_0^1|f(x)|^2x^{\frac{3}{2}}dx, \\
&\|f\|_{[1]}^2=\int_0^1|f(x)|^2(1-x)^{\frac{N}{2}-1}dx.
\end{align*}
Here 
$$u^{[0]}(x)=\omega(x)u(x),\quad u^{[1]}(x)=(1-\omega(x))u(x) $$
with a cut-off function $\omega \in C^{\infty}(\mathbb{R}, [0,1])$ such that 
$\omega(x)=1$ for $0<x \ll 1$ and $\omega(x)=0$ for $0<1-x\ll1$. 
\end{Definition}

Note that 
$$R(t, r_+)=r_+(1+\varepsilon\sin(\sqrt{\lambda}t+\Theta_0)+O(\varepsilon^2)),$$
provided that $\psi$ has been normalized as $\psi(x=1)=1$, and that
the density distribution enjoys the `physical vacuum boundary' condition:
$$\rho(t,r)=
\begin{cases}
C(t)(r_+-r)^{\frac{1}{\gamma-1}}(1+O(r_+-r)) & (0\leq r<r_+) \\
0 & (r_+\leq r)
\end{cases}
$$
with a smooth function $C(t)$ of $t$ such that
$$C(t)=\Big(\frac{\gamma-1}{A\gamma}\frac{Q_+}{r_+^2\kappa_+}\Big)^{\frac{1}{\gamma-1}}+O(\varepsilon).
$$\\

Also we can consider the Cauchy problem
\begin{align*}
&\frac{\partial y}{\partial t}-Jv=0,\qquad
\frac{\partial v}{\partial t}+H_1\mathcal{L}y+H_2=0, \\
&y\Big|_{t=0}=\psi_0(x),\qquad
v\Big|_{t=0}=\psi_1(x).
\end{align*}
Then we have
\begin{Theorem}
Given $T>0$, there exits a small positive $\delta$ such that if
$\psi_0,\psi_1 \in C^{\infty}([0,1])$ satisfy
$$\max_{k\leq\mathfrak{K}}\Big\{\Big\|\Big(\frac{d}{dx}\Big)^k\psi_0\Big\|_{L^{\infty}},
\Big\|\Big(\frac{d}{dx}\Big)^k\psi_1\Big\|_{L^{\infty}}\Big\}\leq \delta, $$
then there exists a unique solution $(y,v)$ of the Cauchy problem
in $C^1([0,T]\times[0,1])$. Here $\mathfrak{K}$ is sufficiently large number.
\end{Theorem}

\section{Metric in the exterior domain}

Let us consider the moving solutions constructed in the preceding
section, which are defined on $0\leq t\leq T, 0<r\leq r_+$.
We discuss on the extension of the metric
onto the exterior vacuum region $r>r_+$.

Keeping in mind the Birkhoff's theorem, we try to patch the Schwarzschild-de Sitter metric
$$ ds^2=\kappa^{\sharp}c^2(dt^{\sharp})^2-
\frac{1}{\kappa^{\sharp}}(dR^{\sharp})^2-(R^{\sharp})^2(d\theta^2+
\sin^2\theta d\phi^2)$$
from the exterior region. Here
$t^{\sharp}=t^{\sharp}(t,r), R^{\sharp}=R^{\sharp}(t,r)$ are smooth functions of $0\leq t\leq T, r_+\leq r\leq r_++\delta$, $\delta$ being a small positive number, and
$$\kappa^{\sharp}=1-\frac{2Gm_+}{c^2R^{\sharp}}-\frac{\Lambda}{3}(R^{\sharp})^2. $$
The patched metric is
$$ds^2=
g_{00}c^2dt^2+2g_{01}cdtdr+g_{11}dr^2+
g_{22}(d\theta^2+\sin^2\theta d\phi^2),
$$
where
\begin{align*}
g_{00}&=\begin{cases}
e^{2F}=\kappa_+e^{-2u/c^2} &\quad (r\leq r_+) \\
\displaystyle \kappa^{\sharp}(\partial_tt^{\sharp})^2-\frac{1}{c^2\kappa^{\sharp}}(\partial_tR^{\sharp})^2 &\quad (r_+<r) 
\end{cases}\\
g_{01}&=
\begin{cases}
0 &\quad (r\leq r_+) \\
\displaystyle c\kappa^{\sharp}(\partial_tt^{\sharp})(\partial_rt^{\sharp})-
\frac{1}{c\kappa^{\sharp}}(\partial_tR^{\sharp})(\partial_rR^{\sharp}) &\quad (r_+<r)
\end{cases}\\
g_{11}&=
\begin{cases}
-e^{2H}=\displaystyle -\Big(1+\frac{V^2}{c^2}-
\frac{2Gm}{c^2R}-\frac{\Lambda}{3}R^2\Big)^{-1}(\partial_rR)^2 &\quad(r\leq r_+) \\
\displaystyle c^2\kappa^{\sharp}(\partial_rt^{\sharp})^2-\frac{1}{\kappa^{\sharp}}(\partial_rR^{\sharp})^2 &\quad (r_+<r)
\end{cases}\\
g_{22}&=\begin{cases}
-R^2 &\quad (r\leq r_+) \\
-(R^{\sharp})^2 &\quad (r_+<r).
\end{cases}
\end{align*}

We require that $R=R^{\sharp}$ and $\partial_rR=\partial_rR^{\sharp}$ along $r=r_+$.
It is necessary for that $g_{22}$ is of class $C^1$. Moreover, by the same way as \cite[Supplementary Remark 4]{TOVdS}, we see that
$$\frac{\partial t^{\sharp}}{\partial t},\quad
\frac{\partial t^{\sharp}}{\partial r}, \quad
\frac{\partial^2t^{\sharp}}{\partial r^2}, \quad
\frac{\partial^2R^{\sharp}}{\partial r^2}$$
at $r=r_++0$ are uniquely determined so that $g_{\mu\nu}$ are of
class $C^1$ across $r=r_+$. By a tedious calculation we have
$$\frac{\partial^2R^{\sharp}}{\partial r^2}\Big|_{r_++0}
-\frac{\partial^2R}{\partial r^2}\Big|_{r_+-0}=\mathcal{A}
\Big(\frac{\partial R}{\partial r}\Big)^2,$$
where
$$\mathcal{A}=-\frac{V^2}{c^2}\Big[
\Big(\frac{Gm_+}{c^2R^2}-\frac{\Lambda}{3}R+
\frac{1}{\sqrt{\kappa_+}}\frac{1}{c^2}\frac{\partial V}{\partial t}\Big)
\Big(1+
\frac{V^2}{c^2}-\frac{2Gm_+}{c^2R}-
\frac{\Lambda}{3}R^2\Big)^{-2}\Big]_{r=r_+-0}.
$$
Since
$$\Big[\Big(\frac{Gm_+}{c^2R^2}-\frac{\Lambda}{3}R+
\frac{1}{\sqrt{\kappa_+}}\frac{1}{c^2}\frac{\partial V}{\partial t}\Big)
\Big(1+
\frac{V^2}{c^2}-\frac{2Gm_+}{c^2R}-
\frac{\Lambda}{3}R^2\Big)^{-2}\Big]_{r=r_+-0}
$$
is near to $\displaystyle \frac{Q_+}{c^2r_+^2\kappa_+^2}\not=0 $,
we see that $\partial^2R^{\sharp}/\partial r^2 \equiv \partial^2R/\partial r^2$
if and only if $V\equiv 0$ at $r=r_+$, which
is the case if the solution under consideration is an equilibrium.

\vspace{20mm}

{\large\bf Acknowledgment} \\

This work is supported by JSPS KAKENHI Grant number JP18K03371.

\end{document}